\documentclass[reqno]{amsart}
\usepackage{hyperref}

\title[Inverse Sturm-Liouville problems]
{Inverse Sturm-Liouville problems with summable potential}

\author[Yu. A. Ashrafyan,   T. N. Harutyunyan]
{Yuri A. Ashrafyan,  Tigran N. Harutyunyan}


\begin{document}

\begin{abstract}
We describe the necessary and sufficient conditions for two sequences $\{\mu_n \} _{n=0}^\infty$ and $\{ a_n \} _{n=0}^\infty$ to be correspondingly the set of eigenvalues and the set of norming constants of a Sturm-Liouville problem with real summable potential $q$ and in advance fixed separated boundary conditions.
\end{abstract}

\maketitle
\numberwithin{equation}{section}
\newtheorem{theorem}{Theorem}[section]
\newtheorem{lemma}[theorem]{Lemma}
\newtheorem{definition}[theorem]{Definition}
\newtheorem{remark}[theorem]{Remark}
\allowdisplaybreaks

\section{Introduction and statements of the results}\label{sec1}

 Let us denote by $L(q, \alpha, \beta)$ the Sturm-Liouville boundary-value problem
\begin{gather}
\ell y\equiv -y'' + q(x)y=\mu y,\quad x\in (0, \pi),\; \mu\in \mathbb{C},\label{eq1.1}\\
y(0) \cos \alpha + y'(0) \sin \alpha = 0, \quad \alpha\in (0, \pi],\label{eq1.2}\\
y(\pi) \cos \beta + y'(\pi) \sin \beta = 0,\quad \beta\in[0, \pi),\label{eq1.3}
\end{gather}
where $q$ is a real-valued, summable function on $[0, \pi]$ (we write $q\in L^1_{\mathbb{R}}[0, \pi]$).
By $L(q, \alpha, \beta)$ we also denote the self-adjoint operator, generated by problem \eqref{eq1.1}-\eqref{eq1.3} (see \cite{Naimark:1969}).
It is well-known, that under these conditions the spectra of the operator  $L(q, \alpha, \beta)$ is discrete and consists of real, simple eigenvalues (see, e.g. \cite{Naimark:1969, Marchenko:1977, Levitan-Sargsyan:1988}), which we denote by $\mu_n=\mu_n(q,\alpha,\beta)=\lambda_n^2(q, \alpha, \beta)$, $n=0, 1, 2, \ldots$, emphasizing the dependence of $\mu_n$ on $q$, $\alpha$ and $\beta$.
We assume that eigenvalues are enumerated in the increasing order, i.e.,
\begin{equation*}
  \mu_0(q, \alpha, \beta) < \mu_1(q, \alpha, \beta) < \dots < \mu_n(q, \alpha, \beta) < \dots.
\end{equation*}

In this article we consider the case $\alpha, \beta \in (0, \pi)$.
It is connected with the circumstance, that in this case the principle term of asymptotics of $\lambda_n = \sqrt{\mu_n}$ is $n$ and the principle term of asymptotics of norming constants $a_n$ (see below \eqref{eq1.4} and \eqref{eq1.7a},\eqref{eq1.7b}) is $\cfrac{\pi}{2}$.
The other three cases: 1) $\alpha = \pi, \beta \in (0,\pi)$, \ 2) $\alpha \in (0, \pi), \beta = 0$, \ 3) $\alpha = \pi, \beta = 0$, need a separate investigation and we do not concern it here.

Let $\varphi(x,\mu) = \varphi(x,\mu,\alpha,q)$ and $\psi(x,\mu) = \psi(x,\mu,\beta,q)$ are the solutions of the equation \eqref{eq1.1}, which satisfy the initial conditions
\begin{gather*}
\varphi(0,\mu,\alpha,q)=1,\qquad \varphi'(0,\mu,\alpha,q)=-\cot\alpha,\\
\psi(\pi,\mu,\beta,q)=1, \qquad  \psi'(\pi,\mu,\beta,q)=-\cot\beta,
\end{gather*}
respectively.
The eigenvalues $\mu_n=\mu_n(q, \alpha, \beta)$, $n=0,1, 2, \ldots$, of $L(q, \alpha, \beta)$ are the zeroes of the characteristic function
\begin{equation*}
  \Delta(\mu): = \varphi(\pi, \mu,\alpha,q)\cot \beta+\varphi'(\pi, \mu,\alpha,q) = -\left( \psi(0, \mu,\beta,q)\cot\alpha+\psi'(0, \mu,\beta,q) \right).
\end{equation*}

It is easy to see that functions $\varphi_n(x):=\varphi(x, \mu_n, \alpha, q)$ and $\psi_n(x):=\psi(x,\mu_n, \beta, q)$, $n=0, 1, 2, \ldots$, are the eigenfunctions, corresponding to the eigenvalue $\mu_n$.
The squares of the $L^2$-norm of these eigenfunctions:
\begin{gather}
a_n=a_n(q,\alpha,\beta):=\int_0^{\pi} |\varphi_n(x)|^2 dx,\quad    n=0, 1, 2, \ldots, \label{eq1.4}\\
b_n=b_n(q,\alpha,\beta):=\int_0^{\pi} |\psi_n(x)|^2 dx,\quad   n=0, 1, 2, \ldots , \label{eq1.5}
\end{gather}
are called norming constants.
The sequences $\{\mu_n \} _{n=0}^\infty$, $\{a_n \} _{n=0}^\infty$ and $\{b_n \} _{n=0}^\infty$ are called spectral data.
The famous theorem of Marchenko (see \cite{Marchenko:1950, Marchenko:1952}) asserts that two sequences $\{\mu_n \} _{n=0}^\infty$ and $\{a_n \} _{n=0}^\infty$ ( or $\{\mu_n \} _{n=0}^\infty$ and $\{b_n \} _{n=0}^\infty$) uniquely determine the problem $L(q, \alpha, \beta)$
\footnote{Recently Ashrafyan has found a new kind of extension of Marchenko theorem, see \cite{Ashrafyan:2017}.}.

In this article we state the question:
\begin{quote}
What kind the sequences $\{\mu_n \} _{n=0}^\infty$ and $\{a_n \} _{n=0}^\infty$ should be, to be the spectral data for a problem $L(q,\alpha,\beta)$ with a $q\in L^1_{\mathbb{R}}[0, \pi]$ and in advance fixed $\alpha$ and $\beta$ from $(0,\pi)$?
\end{quote}

Such a question (but without the condition of fixed $\alpha$ and $\beta$ and for different class of potential $q$ instead of our $q\in L^1_{\mathbb{R}}[0, \pi]$) was considered first by Gelfand and Levitan in work \cite{Gelfand-Levitan:1951} and after in many papers (we refer only some of them: \cite{Gasimov-Levitan:1964, Zhikov:1967, Isaacson-Trubowitz:1983}) and this problem called the inverse Sturm-Liouville problem by "spectral function" (see also, e.g. \cite{Levitan:1984, Freiling-Yurko:2001}).

Our answer to above question is in the following assertion.

\begin{theorem}\label{thm1.1}
 For a real increasing sequence $\{\lambda_n^2 \} _{n=0}^\infty$ and a positive sequence $\{a_n \} _{n=0}^\infty$ to be spectral data for boundary-value problem $L(q, \alpha, \beta)$ with a $q\in L^1_{\mathbb{R}}[0, \pi]$ and fixed $\alpha,\beta \in (0,\pi)$ it is necessary and sufficient that the following relations hold:

\begin{subequations}\label{eq1.6}
1) the sequence $\{\lambda_n \} _{n=0}^\infty$ has asymptotic form
\begin{equation}\label{eq1.6a}
  \lambda_n = n + \cfrac{\omega}{n} + l_n,
\end{equation}
where $\omega = const$,
\begin{equation}\label{eq1.6b}
     l_n = o \left( \cfrac{1}{n} \right), \quad \text{when} \quad n \rightarrow \infty,
\end{equation}
and the function $l(\cdot)$, defined by formula
\begin{equation}\label{eq1.6c}
  l(x) = \sum_{n=1}^{\infty} l_n \sin n x,
\end{equation}
is absolutely continuous on arbitrary segment $[a, b] \subset (0, 2 \pi)$, i.e.
\begin{equation}\label{eq1.6d}
  l \in AC (0, 2 \pi);
\end{equation}
\end{subequations}

\begin{subequations}\label{eq1.7}
2) the sequence $\{a_n \} _{n=0}^\infty$ has asymptotic form
\begin{equation}\label{eq1.7a}
  a_n = \cfrac{\pi}{2} + s_n,
\end{equation}
where
\begin{equation}\label{eq1.7b}
  s_n = o \left( \cfrac{1}{n} \right), \quad \text{when} \quad n \rightarrow \infty,
\end{equation}
and the function $s(\cdot)$, defined by formula
\begin{equation}\label{eq1.7c}
  s(x) = \sum_{n=1}^{\infty} s_n \cos n x,
\end{equation}
is absolutely continuous on arbitrary segment $[a, b] \subset (0, 2 \pi)$, i.e.
\begin{equation}\label{eq1.7d}
  s \in AC (0, 2 \pi);
\end{equation}
\end{subequations}

3)
\begin{equation}\label{eq1.8}
  \cfrac{1}{a_0}-\cfrac{1}{\pi}+\sum_{n=1}^\infty\Big(\cfrac{1}{a_n} - \cfrac{2}{\pi}\Big)=\cot\alpha, \qquad \qquad \qquad \qquad
\end{equation}

4)
\begin{equation}\label{eq1.9}
  \cfrac{a_0}{\big(\pi \prod_{k=1}^\infty  \frac{\mu_k - \mu_0}{k^2}\big)^2}-\cfrac{1}{\pi} + \sum_{n=1}^\infty \left(\cfrac{a_n n^4}{\big( \pi[\mu_0-\mu_n]
 \prod_{k=1, k \neq n}^\infty \frac{\mu_k - \mu_n}{k^2}\big)^2} - \cfrac{2}{\pi}\right)=-\cot\beta.
\end{equation}
\end{theorem}

\vspace{5mm}

In what follows, under condition \eqref{eq1.6} we understand the conditions \eqref{eq1.6a}--\eqref{eq1.6d} and under condition \eqref{eq1.7} the conditions \eqref{eq1.7a}--\eqref{eq1.7d}.

To prove Theorem \ref{thm1.1} we use the following assertion, which has independent interest.

\begin{theorem}\label{thm1.2}
Let $q\in L^1_{\mathbb{R}}[0, \pi]$ and $\alpha,\beta \in (0,\pi)$.
Then for norming constants $a_n=a_n(q, \alpha,\beta)$ and $b_n=b_n(q,\alpha,\beta)$ the following relations are valid
\begin{gather}
\cfrac{1}{{a}_0}-\cfrac{1}{\pi}+\sum_{n=1}^\infty \Big(\cfrac{1}{{a}_n}-\cfrac{2}{\pi}\Big)=\cot\alpha, \label{eq1.10}\\
\cfrac{1}{{b}_0}-\cfrac{1}{\pi}+\sum_{n=1}^\infty \Big(\cfrac{1}{{b}_n}-\cfrac{2}{\pi}\Big)=-\cot\beta. \label{eq1.11}
\end{gather}
\end{theorem}

\vspace{5mm}

Similar results we have obtained in \cite{Ashrafyan-Harutyunyan:2015} for the case $q\in L^2_{\mathbb{R}}[0, \pi]$.
In that case instead of \eqref{eq1.6} we had
\begin{equation}\label{eq1.12}
  \lambda_n = n + \cfrac{\omega}{n} + l_n, \quad \text{where} \quad l_n = \cfrac{\omega_n}{n}, \qquad \{\omega_n\}_{n=0}^\infty \in l^2,
\end{equation}
and instead of \eqref{eq1.7} we had
\begin{equation}\label{eq1.13}
  a_n = \cfrac{\pi}{2} + s_n, \quad \text{where} \quad s_n = \cfrac{\kappa_n}{n}, \qquad \{\kappa_n\}_{n=0}^\infty \in l^2.
\end{equation}

The aim of this paper is to show that when we change the condition $q\in L^2_{\mathbb{R}}[0, \pi]$ to $q\in L^1_{\mathbb{R}}[0, \pi]$, we must change \eqref{eq1.12} by \eqref{eq1.6} and \eqref{eq1.13} by \eqref{eq1.7}.
We should say, that the asymptotics \eqref{eq1.6} and \eqref{eq1.7} have the roots in paper of Zhikov \cite{Zhikov:1967}.
Also we must note that conditions \eqref{eq1.9} and \eqref{eq1.11} are equivalent.
It is a corollary of the fact, that norming constants $b_n = b_n(q, \alpha, \beta), \ n=0, 1, 2, \ldots$, (see \eqref{eq1.5}) can be represented by spectral data $\{\mu_n \} _{n=0}^\infty$ and norming constants $\{a_n \} _{n=0}^\infty$ by the formulae (see \cite{Ashrafyan-Harutyunyan:2015})
\begin{gather}
\frac{1}{b_0} = \frac{a_0}{\pi^2 \big( \prod_{k=1}^\infty \frac{\mu_k - \mu_0}{k^2}\big)^2}, \label{eq1.14}\\
\frac{1}{b_n} = \frac{a_n n^4}{\pi^2 [\mu_0-\mu_n]^2  \big( \prod_{k=1, k \neq n}^\infty \frac{\mu_k - \mu_n}{k^2}\big)^2}, \quad n = 1, 2, \ldots . \label{eq1.15}
\end{gather}
In the same time Theorem \ref{thm1.2} stay valid if we change the case $q\in L^2_{\mathbb{R}}[0, \pi]$ to the case $q\in L^1_{\mathbb{R}}[0, \pi]$

The inverse Sturm-Liouville problem with fixed $\alpha$ and $\beta$ from $(0, \pi)$ and $q\in L^2_{\mathbb{R}}[0, \pi]$ was investigated in \cite{Isaacson-McKean-Trubowitz:1984}, and for the case $\alpha=\pi, \ \beta \in (0, \pi)$ in \cite{Dahlberg-Trubowitz:1984}, in original statement of a question, and the solution of these problems the authors reduced to analysis of inverse Sturm-Liouville problem with Dirichlet boundary conditions ($y(0)=0, \ y(\pi)=0$, which corresponds to $\alpha = \pi, \ \beta = 0$), which was in detail investigated in book \cite{Poshel-Trubowitz:1987}.
We can say, that the problems, solved in \cite{Isaacson-McKean-Trubowitz:1984,  Dahlberg-Trubowitz:1984}, are not the inverse problems by "spectral function", but these problems are deeply connected with the inverse problems by "spectral function", and \cite{Isaacson-McKean-Trubowitz:1984,  Dahlberg-Trubowitz:1984} play an important role in development of inverse problems.
After, the authors of paper \cite{Korotyaev-Chelkak:2009} was studied the solution of inverse Sturm-Liouville problem $L (q, \pi, \beta)$ (i.e. $\alpha = \pi, \beta \in (0, \pi)$) with $q\in L^2_{\mathbb{R}}[0, \pi]$ in terms of eigenvalues and "norming constants" $\{\nu_n \} _{n=0}^\infty$ (see (1.3) in \cite{Korotyaev-Chelkak:2009}), which they introduced for this case
\footnote{We should say, that this "norming constants" had been introduced before in paper \cite{Harutunyan-Nersesyan:2004}, and corresponding uniqueness theorem (see theorem 2.3 in \cite{Korotyaev-Chelkak:2009}) had been proved in \cite{Harutunyan-Nersesyan:2004} (see theorem 3).}.
They proved that the condition \eqref{eq1.11} is necessary for norming constants $b_n(q, \pi, \beta)$ (see \eqref{eq1.5}).
They also proved that the conditions \eqref{eq1.10} and \eqref{eq1.11} are necesssary for norming constants of problem $L(q, \alpha, \beta)$, if $\sin \alpha \neq 0$ and $\sin \beta \neq 0$.
Really they formulated these relations in terms of "norming constants" $\{\nu_n \} _{n=0}^\infty$, but it is easy to verify that these formulations are equivalent.
Thus, Theorem \ref{thm1.2} was proved, for the case $q\in L^2_{\mathbb{R}}[0, \pi]$, with different methods, in \cite{Korotyaev-Chelkak:2009} and \cite{Ashrafyan-Harutyunyan:2015}.
It is easy to see that these proofs remain the same for the case $q\in L^1_{\mathbb{R}}[0, \pi]$.
It is also must be noted, that the relations \eqref{eq1.10} and \eqref{eq1.11} come from the paper \cite{Jodeit-Levitan:1997} of Jodeit and Levitan.

Note that asymptotic behavior of $\{\mu_n\}_{n=0}^\infty$ and  $\{a_n\}_{n=0}^\infty$ are standard conditions for the solvability of the inverse problem (see, e.g., \cite{Gelfand-Levitan:1951, Gasimov-Levitan:1964, Zhikov:1967, Freiling-Yurko:2001}).
The conditions \eqref{eq1.8} and \eqref{eq1.9}, which we add to the conditions \eqref{eq1.6} and \eqref{eq1.7}, guarantee that $\alpha$ and $\beta$, which we construct
during the solution of the inverse problem, are the same that we fixed in advance.
At the same time Theorem \ref{thm1.2} says that the conditions \eqref{eq1.8} and \eqref{eq1.9}, which equivalent to \eqref{eq1.10} and \eqref{eq1.11} are necessary.

\vspace{15mm}

\section{Auxiliary results}\label{sec2}

Consider the function $a(x)$, defined as
\begin{equation}\label{eq2.1}
  a(x) = \sum_{n=0}^{\infty} \left( \cfrac{\cos \lambda_n x}{a_n} - \cfrac{\cos n x}{a_n^0} \right),
\end{equation}
where $a_0^0 = \pi$, $a_n^0 = \cfrac{\pi}{2}$, for $n = 1, 2, \ldots$.

\begin{lemma}\label{lem2.1}
Let the sequences $\{\lambda_n \} _{n=0}^\infty$ and $\{a_n \} _{n=0}^\infty$ have the properties \eqref{eq1.6} and \eqref{eq1.7} correspondingly.
Then $a \in AC (0, 2 \pi)$.
\end{lemma}

\begin{proof}
  We set
  \begin{equation}\label{eq2.2}
    \rho_n = \lambda_n - n = \cfrac{\omega}{n} + l_n = O \left( \cfrac{1}{n} \right).
  \end{equation}
  The general term of the sum \eqref{eq2.1} can be rewritten as follows
  \begin{multline}\label{eq2.3}
    \cfrac{\cos \lambda_n x}{a_n} - \cfrac{\cos n x}{a_n^0} =
    \cfrac{\cos \lambda_n x}{a_n} - \cfrac{\cos n x}{a_n} + \cfrac{\cos n x}{a_n} - \cfrac{\cos n x}{a_n^0} = \\
    = \cfrac{1}{a_n} \left( \cos \lambda_n x - \cos n x \right) + \left( \cfrac{1}{a_n} - \cfrac{1}{a_n^0} \right) \cos n x.
  \end{multline}
  So we can rewrite the series \eqref{eq2.1} in the following form
  \begin{equation}\label{eq2.10}
    a(x) = \sum_{n=0}^{\infty} \cfrac{1}{a_n} \left( \cos \lambda_n x - \cos n x \right) +
            \sum_{n=0}^{\infty} \left( \cfrac{1}{a_n} - \cfrac{1}{a_n^0} \right) \cos n x.
  \end{equation}
  The difference $\left( \cos \lambda_n x - \cos n x \right)$ can be represented as follows
  \begin{multline}\label{eq2.4}
    \cos \lambda_n x - \cos n x =
    \cos (n + \rho_n) x - \cos n x = \\
    = \cos n x \cos \rho_n x - \sin n x \sin \rho_n x - \cos n x = \\
    = - \cos n x (1 - \cos \rho_n x) - \sin n x \sin \rho_n x = \\
    = - 2 \sin^2 \cfrac{\rho_n x}{2} \cos n x - \rho_n x \sin n x - (\sin \rho_n x - \rho_n x) \sin n x
  \end{multline}
  and for the difference $ \left( \cfrac{1}{a_n} - \cfrac{1}{a_n^0} \right)$ we have
  \begin{equation}\label{eq2.5}
    \cfrac{1}{a_n} - \cfrac{1}{a_n^0} =  \cfrac{1}{\frac{\pi}{2} + s_n} - \cfrac{1}{\frac{\pi}{2}} = - \cfrac{2}{\pi} \cdot \cfrac{s_n}{1 + \frac{\pi}{2}s_n}
    = - \cfrac{2}{\pi} \cdot s_n + q_n,
  \end{equation}
  where $q_n = o \left( \cfrac{1}{n^2} \right)$.
  From the latter relation \eqref{eq2.5} it follows, that for sufficiently large $n$ we have
  \begin{equation}\label{eq2.6}
    \cfrac{1}{a_n} = \cfrac{1}{a_n^0} + o \left( \cfrac{1}{n} \right) = \cfrac{2}{\pi} + o \left( \cfrac{1}{n} \right).
  \end{equation}
  And then, according to the relations \eqref{eq2.2} and \eqref{eq2.6}, we obtain
  \begin{multline}\label{eq2.7}
    \cfrac{1}{a_n}(- \rho_n x) \sin n x = \left[ - \cfrac{2}{\pi} + o \left( \cfrac{1}{n} \right) \right] \left[ \cfrac{\omega}{n} + l_n \right] x \sin n x = \\
    = - \cfrac{2}{\pi} \omega x \cfrac{\sin n x}{n} - \cfrac{2}{\pi} l_n x \sin n x + r_n \sin n x,
  \end{multline}
  where $r_n = o \left( \cfrac{1}{n^2} \right)$.

  Since $\rho_n = O \left( \cfrac{1}{n} \right)$ and $\sin y - y = O \left( n^3 \right)$ for $y$ close to zero, then
  \begin{gather}
    \sin \rho_n x - \rho_n x = O \left( \left(\rho_n x\right)^3 \right) = O \left( \cfrac{1}{n^3} \right), \label{eq2.8} \\
    \sin^2 \cfrac{\rho_n x}{2} = O \left( \left(\cfrac{\rho_n x}{2}\right)^2 \right) = O \left( \cfrac{1}{n^2} \right). \label{eq2.9}
  \end{gather}

  Thus, taking into account the relations \eqref{eq2.4}--\eqref{eq2.9} we can rewrite the function \eqref{eq2.10} as follows
  \begin{equation*}
    a(x) = a_1(x) + a_2(x),
  \end{equation*}
  where
  \begin{gather*}
    a_1(x) =  - \cfrac{2 \omega x}{\pi} \sum_{n=1}^{\infty} \cfrac{\sin n x}{n} - \cfrac{2 x}{\pi} \sum_{n=1}^{\infty} l_n \sin n x -
    \cfrac{2}{\pi} \sum_{n=1}^{\infty}  s_n \cos n x, \\
    a_2(x) = - \sum_{n=1}^{\infty} \cfrac{1}{a_n} \left(\sin \rho_n x - \rho_n x \right) \sin n x - \sum_{n=0}^{\infty} \cfrac{1}{a_n} \sin^2 \cfrac{\rho_n x}{2} \cos n x + \\
    + \sum_{n=0}^{\infty} q_n \cos n x + \sum_{n=1}^{\infty} r_n \sin n x.
  \end{gather*}

  Since the first sum $\sum_{n=1}^{\infty} \cfrac{\sin n x}{n} = \cfrac{\pi - x}{2}$, for $x \in (0, 2 \pi)$, then, in particular, it belongs to $AC(0, 2 \pi)$,
  and the second sum $\sum_{n=1}^{\infty} l_n \sin n x$ also belongs to $AC(0, 2 \pi)$ according to the condition \eqref{eq1.6d} of the lemma.
  In its turn the sum $\sum_{n=1}^{\infty} s_n \cos n x$ belongs to $AC(0, 2 \pi)$ according to the condition \eqref{eq1.7d} of the lemma.
  The other four sums of $a_2(x)$ converge absolutely and uniformly on $[0, 2 \pi]$ and are continuous differentiable functions, and hence belong to $AC(0, 2 \pi)$.

  These complete the proof.
\end{proof}

\begin{lemma}\label{lem2.2}
Let the sequences $\{\lambda_n \} _{n=0}^\infty$ and $\{a_n \} _{n=0}^\infty$ have the properties \eqref{eq1.6} and \eqref{eq1.7}, \eqref{eq1.8} correspondingly.
Then the function $F$, defined in triangle $0 \leq t \leq x \leq \pi$ by formula
\begin{equation}\label{eq2.11}
  F(x, t) = \sum_{n=0}^{\infty} \left( \cfrac{\cos \lambda_n x \cos \lambda_n t}{a_n} - \cfrac{\cos n x \cos n t}{a_n^0} \right),
\end{equation}
is absolutely continuous function with respect to each variable and function
\begin{equation*}
  f(x) := \cfrac{d}{dx} F(x, x),
\end{equation*}
is summable on $(0, \pi)$, i.e. $f \in L^1(0, \pi)$.
\end{lemma}

\begin{proof}
It is easy to see that
\begin{equation}\label{eq2.12}
  F(x, t) = \cfrac{1}{2} \left[ a(x + t) + a(x - t) \right].
\end{equation}
Since,  $a \in AC (0, 2 \pi)$, then we can infer that the function $F(x, t)$ with respect to both of the variables has the same smoothness as $a(x)$.
For the functions $F(x, x)$ we have
\begin{equation}\label{eq2.13}
  F(x, x) = \cfrac{1}{2} \left[ a(2 x) + a(0) \right].
\end{equation}
According to \eqref{eq1.8} and \eqref{eq2.1}
\begin{equation*}
  a(0) = \sum_{n=0}^{\infty} \left( \cfrac{1}{a_n} - \cfrac{1}{a_n^0} \right) = \cot \alpha,
\end{equation*}
for $\alpha \in (0, \pi)$.
Hence $a(0)$ has a sense and, therefore, $F(x, x)$ too.
Besides this
\begin{equation}\label{eq2.14}
  \cfrac{d}{dx} F(x, x) = \cfrac{1}{2} \cfrac{d}{dx} a(2 x),
\end{equation}
and since $a \in AC(0, 2\pi)$, then the function $a(2x)$ belongs to $AC(0, \pi)$, and its derivative belongs to $L^1_\mathbb{R}(0, \pi)$, i.e. the function
$\cfrac{d}{dx} F(x, x)$ belongs to $L^1_\mathbb{R}(0, \pi)$.
\end{proof}

\vspace{15mm}

\section{proof of the Theorem \ref{thm1.1}}\label{sec3}

\textbf{The proof of necessity.}
If $\{\lambda_n^2 \} _{n=0}^\infty$ are the eigenvalues and $\{a_n \} _{n=0}^\infty$ are the norming constants of the problem $L(q, \alpha, \beta)$, then for $\mu_n = \lambda_n^2$ the asymptotics \eqref{eq1.6} was proved in \cite{Harutyunyan:2016}, and for $a_n$ the asymptotics \eqref{eq1.7} was proved in \cite{Harutyunyan-Pahlevanyan:2016}.
The necessity of connections \eqref{eq1.8} and \eqref{eq1.9} was proved in \cite{Ashrafyan-Harutyunyan:2015} for $q\in L_{\mathbb{R}}^2[0, \pi]$, and the same proof is true for the case $q\in L_{\mathbb{R}}^1[0, \pi]$.

\vspace{5mm}

\textbf{The proof of sufficiency.}
In \cite{Gasimov-Levitan:1964} there is a proof of such assertion:

\begin{theorem}[\cite{Gasimov-Levitan:1964}]\label{thm3.1}
For real numbers $\{\lambda_n^2 \} _{n=0}^\infty$ and $\{a_n \} _{n=0}^\infty$ to be the spectral data for a certain boundary-value problem $L(q, \alpha, \beta)$ with $q\in L_{\mathbb{R}}^1[0, \pi]$, ($\alpha,\beta \in \left( 0,\pi \right)$), it is necessary and sufficient that  relations \eqref{eq1.6a}--\eqref{eq1.6b} and \eqref{eq1.7a}--\eqref{eq1.7b} hold, and the function $F(\cdot , \cdot)$ has partial derivatives, which are summable with respect to each variable.
\end{theorem}

Thus, if we have a real sequence $\{\mu_n \}_{n=0}^\infty = \{\lambda_n^2 \}_{n=0}^\infty$, which has the asymptotic representation \eqref{eq1.6a}--\eqref{eq1.6b} and a positive sequence $\{a_n \}_{n=0}^\infty$, which has the asymptotic representation \eqref{eq1.7a}--\eqref{eq1.7b}, then, according to the Theorem \ref{thm3.1}, there exist a function $q\in L^1_{\mathbb{R}}[0, \pi]$ and some constants $\tilde{\alpha},\tilde{\beta} \in (0,\pi)$ such that $\lambda_n^2$, $n=0,1,2,\ldots$, are the eigenvalues and ${a}_n$, $n=0,1,2,\ldots$, are norming constants of a Sturm-Liouville problem $L(q,\tilde{\alpha},\tilde{\beta})$.

The function $q(x)$ and constants $\tilde{\alpha}, \tilde{\beta}$ are obtained on the way of solving the inverse problem by Gelfand-Levitan method.
The algorithm of that method is as follows:

First we define the function $F(x,t)$ by formula
\begin{equation}\label{eq3.1}
F(x,t)=\sum_{n=0}^\infty \left( \cfrac{\cos \lambda_n x \cos \lambda_n t }{a_n} - \cfrac{\cos n x \cos n t }{a_n^0} \right).
\end{equation}
Note that this function is defined by $\{\lambda_n \}_{n=0}^\infty$ and $\{{a}_n \}_{n=0}^\infty$ uniquely.
Then we solve Gelfand-Levitan integral equation \cite{Marchenko:1950, Gelfand-Levitan:1951, Gasimov-Levitan:1964, Freiling-Yurko:2001}
\begin{equation}\label{eq3.2}
G(x,t)+F(x,t)+\int_0^x G(x,s)F(s,t)ds=0, \qquad 0 \leq t \leq x,
\end{equation}
where $G(x,\cdot)$ is unknown function.
Find function $G(x,t)$, with the help of which we construct a function
\begin{equation}\label{eq3.3}
\varphi(x, \lambda^2) = \cos \lambda x + \int_0^x G(x,t) \cos \lambda t dt,
\end{equation}
which is defined for all $\lambda \in \mathbb{C}$. It is proved (see \cite{Gasimov-Levitan:1964}) that
\begin{equation}\label{eq3.4}
- {\varphi}''(x,\lambda^2)+\Big( 2 \cfrac{d}{dx}G(x,x) \Big) {\varphi}(x,\lambda^2)
 = \lambda^2 {\varphi}(x,\lambda^2),
\end{equation}
almost everywhere on $(0,\pi)$, and
\begin{gather*}
\varphi(0,\lambda^2)=1,\\
\varphi'(0,\lambda^2)=G(0,0).
\end{gather*}
If we denote
\begin{equation}\label{eq3.5}
G(0,0)=-\cot \tilde \alpha,
\end{equation}
then the solution \eqref{eq3.3} of  equation \eqref{eq3.4} will satisfy the boundary condition \eqref{eq1.2}
\begin{equation*}
  \varphi(0,\lambda^2) \cot \tilde \alpha + \varphi'(0,\lambda^2) =0
\end{equation*}
for all $\lambda \in \mathbb{C}$.
Since from \eqref{eq3.2} follows that $G(0,0)=-F(0,0)$ and from \eqref{eq3.1} follows that $F(0,0)= \displaystyle \sum_{n=0}^\infty \left( \cfrac{1}{a_n}-\cfrac{1}{a_n^0}\right)$, hence we get
\begin{equation}\label{eq3.6}
  \sum_{n=0}^\infty \left( \cfrac{1}{a_n}-\cfrac{1}{a_n^0}\right)=\cot \tilde \alpha.
\end{equation}
From the relation \eqref{eq3.6} and our condition \eqref{eq1.8} on the sequence $\{a_n \}_{n=0}^\infty$ we find that $\tilde \alpha = \alpha$.

It is also proved (see, e.g., \cite{Gasimov-Levitan:1964}) that the expression
\begin{equation*}
   \cfrac {\varphi'_n (\pi)}{\varphi_n (\pi)} = \cfrac {\varphi'(\pi, \lambda_n^2)}{\varphi(\pi,\lambda_n^2)}
\end{equation*}
is a constant (i.e. does not depend on $n$), which we will denote by $-\cot \tilde{\beta}$.
Thus the functions $\varphi(x,\lambda_n^2), n=0,1,2,\ldots$, are the eigenfunctions of a problem $L(q,\tilde{\alpha},\tilde{\beta})$, where $q(x) = 2 \cfrac{d}{dx} G(x,x)$, $\tilde{\alpha}$ is in advance given $\alpha$ and we should have $\tilde{\beta}$ equals $\beta$.
We know from the Theorem \ref{thm1.2}, that for problem $L(q,\alpha,\tilde{\beta})$ it holds
\begin{equation*}
  \cfrac{1}{b_0}-\cfrac{1}{\pi}+\sum_{n=1}^\infty \left( \cfrac{1}{b_n}- \cfrac{2}{\pi} \right) = - \cot \tilde{\beta}.
\end{equation*}
Thus, if we obtain condition \eqref{eq1.11}, then we guarantee that $\tilde{\beta}=\beta$.
But \eqref{eq1.11} deals with the norming constants $b_n$, which are not independent.
We have shown in the paper \cite{Ashrafyan-Harutyunyan:2015} that we can represent $b_n$ by $a_n$ and $\{\mu_k \}_{k=0}^\infty$ by formulae \eqref{eq1.14}, \eqref{eq1.15}.
Therefore, instead of \eqref{eq1.11}, we obtain the condition in the form \eqref{eq1.9}.

This completes the proof.

\vspace{15mm}

\textbf{Acknowledgment}

This work was supported by the RA MES State Committee of  Science, in the frames of the research project No.15T-1A392.

\end{document}